\documentclass{amsart}
\usepackage{amssymb} 
\usepackage{amsmath} 
\usepackage{comment}
\usepackage[OT2,T1]{fontenc}
\DeclareSymbolFont{cyrletters}{OT2}{wncyr}{m}{n}
\DeclareMathSymbol{\Sha}{\mathalpha}{cyrletters}{"58}
\DeclareMathOperator{\spn}{Span}
\DeclareMathOperator{\GL}{GL}

\DeclareMathOperator{\ord}{ord}

\newcommand{\Q}{\mathbb{Q}}

\newcommand{\R}{\mathbb{R}}
\newcommand{\F}{\mathbb{F}}
\newcommand{\PP}{\mathbb{P}}

\newcommand{\cF}{\mathcal{F}}
\newcommand{\cB}{\mathcal{B}}

\newcommand{\cT}{\mathcal{T}}

\newcommand{\cC}{\mathcal C}
\newcommand{\cU}{\mathcal{U}}
\newcommand{\cV}{\mathcal{V}}
\newcommand{\cW}{\mathcal{W}}

\newcommand{\vv}{\mathbf{v}}

\newcommand{\xx}{\mathbf{x}}
\newcommand{\yy}{\mathbf{y}}
\newcommand{\ff}{\mathbf{f}}
\newcommand{\fg}{\mathbf{g}}
\newcommand{\bz}{\boldsymbol{\zeta}}
\newcommand{\bx}{\boldsymbol{\xi}}

\begin {document}

\newtheorem{thm}{Theorem}

\newtheorem{lem}{Lemma}[section]
\newtheorem{prop}[lem]{Proposition}

\newtheorem{cor}[lem]{Corollary}

\theoremstyle{definition}

\theoremstyle{remark}

\title[Mordell--Weil in High Dimension]{
A Mordell--Weil Theorem for Cubic Hypersurfaces
of High Dimension
}
\author{Stefanos Papanikolopoulos} 
\author{Samir Siksek}
\address{Mathematics Institute\\
	University of Warwick\\
	Coventry\\
	United Kingdom}

\email{s.siksek@warwick.ac.uk}
\date{\today}
\thanks{
Siksek is supported by the 
EPSRC \emph{LMF: L-Functions and Modular Forms} Programme Grant
EP/K034383/1.
}

\keywords{cubic hypersurfaces, rational points, Mordell--Weil problem}
\subjclass[2010]{Primary 14G35, Secondary 14J70}

\begin {abstract}
Let $X/\Q$ be a smooth cubic hypersurface of dimension $n \ge 1$.
It is well-known that new rational points
may be obtained from old ones by secant and tangent constructions.
In view of the Mordell--Weil theorem for $n=1$, 
Manin (1968) asked if there exists a finite set $S$
from which all other rational points can be thus 
obtained.
We give an affirmative
answer for $n \ge 48$, showing in fact that we can take the
generating set $S$
to consist of just one point. Our proof makes use of a weak
approximation theorem 
due to 
Skinner, a theorem of 
Browning, Dietmann and Heath-Brown on the existence
of rational points on the intersection
of a quadric and cubic in large dimension, and some
elementary ideas from differential geometry, 
algebraic geometry and numerical analysis.
\end {abstract}
\maketitle

\section{Introduction}


Let $X \subset \PP^{n+1}$ 
be a smooth cubic hypersurface over $\Q$ of dimension $n$.
Let $\ell \subset \PP^{n+1}$ be a line defined over $\Q$.
If $\ell$ is not contained in $X$ then
$\ell \cdot X=P+Q+R$ where $P$, $Q$, $R\in X$. If
any two of $P$, $Q$, $R$ are rational then so is the third.
If $S \subseteq X(\Q)$, we write $\spn(S)$
for the subset of $X(\Q)$ generated from $S$ by successive
secant and tangent constructions. 
More formally, we define a sequence
\[
S=S_0 \subseteq S_1 \subseteq S_2 \subseteq \cdots \subseteq X(\Q)
\]
by letting $S_{n+1}$ be the set of points $R\in X(\Q)$
such that either $R \in S_n$,
or for some $\Q$-line $\ell \not \subset X$ we have
$\ell \cdot X=P+Q+R$ where $P$, $Q \in S_n$. Then
$\spn(S):= \cup S_n$. Manin \cite[page 3]{Ma1}
asks if there is some finite
subset $S \subset X(\Q)$ such that $\spn(S)=X(\Q)$. 
\begin{thm}\label{thm:main}
Let $X$ be a smooth cubic hypersurface of dimension $n \ge 48$
defined over $\Q$. Then there exists a point $A \in X(\Q)$
such that $\spn(A)=X(\Q)$.
\end{thm}

We are grateful to Tim Browning, 
Simon Rydin Myerson,
Michael Stoll
and Damiano Testa for useful discussions.

\bigskip

\subsection{Notation} Throughout $X \subset \PP^{n+1}$
a smooth cubic hypersurface of dimension $n$ defined over $\Q$
(for now $n \ge 2$).
Thus there is some non-zero homogeneous cubic
polynomial $F \in \Q[x_0,\dotsc,x_{n+1}]$ such that $X$
is given by the equation
\begin{equation}\label{eqn:F}
X \; : \; F(x_0,\dotsc,x_{n+1})=0.
\end{equation}
For $P\in X$ we let $T_P{X}$ denote the tangent
plane to $X$ at $P$:  
\[
T_P{X} \; : \; \nabla{F}(P) \cdot (x_0,\dotsc,x_{n+1})=0.
\]
The Gauss map on $X$ sends $P$ to $T_P{X} \in {\PP^{n+1}}^*$.
We let 
$X_P:=X \cap T_P{X}$. Thus 
\[
X_P \; : \; 
\begin{cases}
F(x_0,\dotsc,x_n)=0,\\
\nabla{F}(P) \cdot (x_0,\dotsc,x_n)=0.
\end{cases}
\]
In Section~\ref{sec:diffgeom} we introduce the second fundamental
form $\Pi_P{X}$, and the Hessian $H_F(P)$.
We write $\mathbb{G}(n+1,1)$ for the Grassmannian parametrizing
lines in $\PP^{n+1}$. Throughout the terms \lq open\rq\ and \lq closed\rq\
will be with respect to the real topology, unless prefixed by \lq Zariski\rq. 

\subsection{A sketch of the proof of Theorem~\ref{thm:main}}
We show in Section~\ref{sec:tangentplane} that if $B \in X(\Q)$ is
not an Eckardt point then 
$X_B(\Q) \subseteq \spn(B)$ (the definition of Eckardt points is given in
Section~\ref{sec:diffgeom}). 
Fix $B \in X(\Q)$ that is non-Eckardt.
Given $D \in X(\Q)$, we ask if there is $C \in X_B(\Q)$
such that $D \in X_C(\Q)$? 
If so, then provided $C$ is non-Eckardt, we have
$D \in \spn(C) \subseteq \spn(B)$. 
The answer to this question is positive provided the variety $Y_{B,D} \subset \PP^{n+1}$ given by
\begin{equation}\label{eqn:intersection}
Y_{B,D} \; : \; 
\begin{cases}
F(x_0,\dotsc,x_{n+1})=0\\
\nabla{F}(x_0,\dotsc,x_{n+1}) \cdot D=0\\
\nabla{F}(B) \cdot (x_0,\dotsc,x_{n+1})=0.
\end{cases}
\end{equation}
has a rational point.
A theorem of Browning, Dietmann and Heath-Brown
allows us to deduce the existence of a rational point under some conditions,
the most important being that $n$ is large, and that $Y_{B,D}$ 
has a smooth real point. By considering the second fundamental form,
and using a theorem on weak approximation for cubic hypersurfaces
due to Skinner, we shall show the existence
of a point $B \in X(\Q)$ and a non-empty open $U \subseteq X(\R)$,
so that $Y_{B,D}$ 
has a smooth real point for all $D \in U$. It follows (with a little care) that
$U \cap X(\Q) \subseteq \spn(B)$. Once the existence of such
a set $U$ is established, we use Mordell--Weil operations to enlarge
$U$ and quickly complete the proof of Theorem~\ref{thm:main}.

\section{Some results from analytic number theory}

\subsection{Weak Approximation}
We shall need the following theorem of Skinner \cite{Skinner}.
\begin{thm}[Skinner]\label{thm:Skinner}
Suppose $n \ge 15$. Then $X$ satisfies weak approximation.
\end{thm}
This means that $X(\Q)$ is dense in $X(A_\Q)$ where $A_\Q$
denotes the adeles.
It follows that $X(\Q)$ is dense in $X(\R)$; a fact we use
repeatedly in the proof of Theorem~\ref{thm:main}.

\begin{cor}\label{cor:tangentbundle}
Suppose $n \ge 15$. Let $U$, $V \subseteq X(\R)$ be disjoint open sets.
Let $A^\prime \in U$, $B^\prime \in V$, and let
$\ell^\prime \not \subset X$ be an $\R$-line 
such that $\ell^\prime \cdot X=2A^\prime+B^\prime$.
Then there are 
$A \in U \cap X(\Q)$, $B \in V \cap X(\Q)$ and a $\Q$-line
$\ell \not \subset X$ such 
that $\ell \cdot X=2A+B$.
\end{cor}
\begin{proof}
The projectivized tangent bundle $\cT_X$ of $X$ parametrizes
pairs $(P,\ell)$ with $P\in X$ and $\ell$ a line
tangent to $X$ at $P$. 
We make use of the fact that $\cT_X$
is locally trivial; thus there is a Zariski open $\cU$ containing $A^\prime$,
and a local isomorphism $\varphi \; : \; \cU \times \PP^{n-1} \rightarrow \cT_X$
such that $\varphi(P,\alpha)=(P,\ell_{P,\alpha})$ where $\ell_{P,\alpha}$
is a line tangent to $X$ at $P$. 
Moreover, as $A^\prime$ is real we
take $\varphi$ to be defined over $\R$. 
Let $W=\cU(\R) \cap U$ which is
necessarily an open neighbourhood of $A^\prime$. Let
$\alpha \in \PP^{n-1}(\R)$ so that $\ell^\prime=\ell_{A^\prime,\alpha}$.
By 
Theorem~\ref{thm:Skinner} we can find 
$\{A_i\} \subset  W \cap X(\Q)$ converging to $A^\prime$.
Write $\ell_i=\ell_{A_i,\alpha}$. Then $\{\ell_i\}$ converges to $\ell^\prime$
in $\mathbb{G}(n+1,1)(\R)$. In particular, for sufficiently large $i$,
the line $\ell_i$ meets $V$. Let $A=A_i \in U \cap X(\Q)$ 
for any such large $i$. Choosing a line $\ell/\Q$ tangent to $X$ at $A$
that sufficiently approximates $\ell_i$
completes the proof.
\end{proof}

\subsection{Intersections of a cubic with a quadric}
Let $Q$, $C \in \Q[x_1,\dots,x_k]$ be a pair of forms
of degrees $2$ and $3$ respectively, such that 
\[
Y \; : \; C(x_1,\dotsc,x_k)=0, \qquad Q(x_1,\dotsc,x_k)=0
\]
is a complete intersection $Y \subset \PP^{k-1}$. 
Using the circle method, Browning,
Dietmann and Heath-Brown establish various sufficient 
conditions for $Y$ to have a rational
point. 
We recount
one of their theorems \cite[Theorem 4]{BDH} which will be essential to our
proof of Theorem~\ref{thm:main}. 
Define $\ord_Q(C)$ as the least non-negative
integer
$m$ such that $C=LQ+C^\prime$, with linear form $L \in \Q[x_1,\dotsc,x_k]$, 
and such that there is a matrix $\mathbf{T} \in \GL_k(\Q)$
with $C^\prime(\mathbf{T}(x_1,\dots,x_k)) \in \Q[x_1,\dotsc,x_m]$.

\begin{thm}[Browning, Dietmann and Heath-Brown]\label{thm:BDH}
With notation as above, 
suppose $k \ge 49$ and $\ord_Q(C) \ge 17$. 
If $Y$ has a smooth real point
then $Y(\Q) \ne \emptyset$.
\end{thm}

\begin{cor}\label{cor:BDH}
Let $f$, $q$, $l \in \Q[x_0,\dotsc,x_{n+1}]$, be forms of
degree $3$, $2$, $1$. Write
\[
Z \; : \; f(x_0,\dotsc,x_{n+1})=q(x_0,\dotsc,x_{n+1})=l(x_0,\dotsc,x_{n+1})=0
\]
for their common locus of zeros in $\PP^{n+1}$. Suppose
that 
\begin{enumerate}
\item[(i)] the cubic hypersurface in $\PP^{n+1}$ defined by $f$
is smooth;
\item[(ii)] $Z$ has a smooth real point;
\item[(iii)] $n \ge 48$.
\end{enumerate}
Then $Z$ has a rational point.
\end{cor}
\begin{proof}
By a non-singular change of variable, we may suppose that
$l=x_0$. Let 
\[
f^\prime(x_1,\dotsc,x_{n+1})=f(0,x_1,\dotsc,x_{n+1}), \qquad
q^\prime(x_1,\dotsc,x_{n+1})=q(0,x_1,\dotsc,x_{n+1}).
\]
We may therefore consider $Z$ as being given
in $\PP^n$ as the common locus of $f^\prime=q^\prime=0$.
Suppose $\ord_{q^\prime}(f^\prime) \le 16$. Then a further non-singular
change of
variables allows us to write
\[
f=x_{0} q_0 + l^\prime q^\prime+ h(x_1,\dots,x_{16}). 
\]
where $q_0$ is a quadratic form, $l^\prime$ is a linear form, 
and $h$ is a cubic form.
Now as $n \ge 48$, there is a common zero in $\mathbb{P}^{n+1}$ to
\[
x_0=x_1=\cdots=x_{16}=l^\prime=q_0=q^\prime=0.
\]
This gives a singular point
on the cubic hypersurface $f=0$ in $\PP^{n+1}$ contradicting
(i). We may thus suppose that $\ord_{q^\prime}(f^\prime) \ge 17$.
A similar argument shows that $f^\prime=q^\prime=0$ defines
a complete intersection in $\PP^{n}$.
By (ii) this intersection has a smooth real point.
Applying Theorem~\ref{thm:BDH} with $k=n+1$ completes the proof.
\end{proof}

\section{A Numerical Stability Criterion}

\subsection{Newton--Raphson}
We need a rigorous version of the multivariate Newton--Raphson
method. The following result is part of Theorem 5.3.2 of 
\cite{numerical}. Here $\lVert \cdot \rVert$ denotes the usual 
Euclidean norm (both for vectors and for
matrices). For differentiable  
$\ff=(f_1,\dots,f_n) \, : \, \R^n
\rightarrow \R^n$, denote the Jacobian matrix by $J_\ff$: 
\[
J_\ff:=\left(\frac{\partial{f_i}}{\partial{x_j}} \right)_{i,j=1,\dotsc,n} .
\]
\begin{thm}\label{thm:nr}
Let $\cC \subseteq \R^n$ be open, $\cC_0$ be convex 
such that $\overline{\cC_0} \subseteq \cC$, and let 
$\ff \, : \, \R^n \rightarrow \R^n$ be differentiable
for all $\xx \in \cC_0$ and continuous for all $\xx \in \cC$.

For $\xx_0 \in \cC_0$ let $r$, $\alpha$, $\beta$, $\gamma$, $h$
be given with the following properties:
\[
\cB_r(\xx_0):=\{ \, \xx \, : \, \lVert \xx-\xx_0 \rVert< r \}
\subseteq \cC_0,  \qquad
h:=\alpha \beta \gamma/2 <1, \qquad r:=\alpha/(1-h),
\]
and let $\ff$ satisfy:
\begin{enumerate}
\item[(i)] $\lVert J_\ff(\xx)-J_\ff(\yy) \rVert \le \gamma \lVert \xx-\yy \rVert$
for all $\xx$, $\yy \in \cC_0$;
\item[(ii)] $J_\ff(\xx)^{-1}$ exists and satisfies 
$\lVert J_\ff(\xx)^{-1} \rVert \le \beta$ for all $\xx \in \cC_0$;
\item[(iii)] $\lVert \ff(\xx_0) \cdot J_\ff(\xx_0)^{-1} \rVert \le \alpha$.
\end{enumerate}
Then beginning at $\xx_0$ each point 
\[
\xx_{k+1}=\xx_k-\ff(\xx_k) \cdot J_\ff(\xx_k)^{-1} \, , \, \qquad k=0,1,2,\dots
\]
is well-defined and belongs to $\cB_r(\xx_0)$. Moreover the limit
$\lim_{k \rightarrow \infty} \xx_k=\bx$ exists,
belongs to $\overline{\cB_r(\xx_0)}$ and 
satisfies $\ff(\boldsymbol{\xi})=\mathbf{0}$.
\end{thm}
\subsection{Stability}
For $f \in \R[x_1,\dotsc,x_n]$ we shall let $\lVert f \rVert$
denote the maximum of the absolute values of the coefficients of $f$.
\begin{lem}\label{lem:nr}
Let $g_1,\dotsc,g_m \in \R[x_1,\dotsc,x_n]$ be polynomials with
$m \le n$. Let $\bz \in \R^n$ 
be a common zero of $g_1,\dots,g_m$,
such that $\nabla{g_1}(\bz),\dotsc,\nabla{g_m}(\bz)$ are linearly
independent. Let $\varepsilon >0$ be given. There is
$\delta>0$ such that if $f_1,\dotsc,f_m \in \R[x_1,\dotsc,x_n]$
satisfy $\lVert f_i-g_i \rVert<\delta$, then there is
$\bx \in \R^n$ such that
\begin{enumerate}
\item[(a)] $\bx$ is a common zero to $f_1,\dotsc,f_m$;
\item[(b)] $\nabla{f_1}(\bx),\dotsc,\nabla{f_m}(\bx)$
are linearly independent;
\item[(c)] $\lVert \bx-\bz \rVert< \varepsilon$.
\end{enumerate}
\end{lem}
\begin{proof}
Choose $\vv_{m+1},\dotsc,\vv_n \in \R^n$ so that
$\nabla{g_1}(\bz),\dotsc,\nabla{g_m}(\bz),\vv_{m+1},\dotsc,\vv_n$
is a basis. Let 
\[
g_i(\xx)=\vv_i \cdot (\xx-\bz), \qquad i=m+1,\dots,n. 
\] 
Then $\bz$ is a common zero to $g_1,\dots,g_n$
and $\nabla{g_1}(\bz),\dotsc,\nabla{g_n}(\bz)$
are linearly independent. Let $\fg=(g_1,\dotsc,g_n)$.
Then $J_\fg(\bz)$ is invertible.
We shall fix 
$f_i=g_i$ for $i=m+1,\dotsc,n$, and let $\ff=(f_1,\dotsc,f_n)$.
We will apply Theorem~\ref{thm:nr} with $\xx_0=\bz$.
There is some $\delta_0>0$ such that if $\lVert f_i - g_i \rVert<\delta_0$
then $J_\ff(\xx_0)$ is invertible. Choose $0<r_0 \le \varepsilon$
so that condition (ii) of the theorem
is satisfied for some $\beta>0$, with $\cC_0=\overline{\cB_{r_0}(\xx_0)}$.
Condition (i) holds for some $\gamma>0$ by the multivariate
Taylor Theorem. Let $\alpha=\lVert \ff(\xx_0) \cdot J_\ff(\xx_0)^{-1} \rVert$,
which depends on $\ff$. Now $\fg(\xx_0)=\mathbf{0}$, so
clearly if $\delta \rightarrow 0$,
then $\alpha \rightarrow 0$. Therefore for sufficiently small
$\delta<\delta_0$, we have $h:=\alpha \beta \gamma/2<1$ and 
$r:=\alpha/(1-h)<r_0$. By the theorem, there
is $\bx \in \overline{\cB_r(\xx_0)}$ such that $\ff(\bx)=\mathbf{0}$.
By construction $\bx$ satisfies (a), (b), (c).
\end{proof}

\subsection{Smooth Real Points on the Varieties $Y_{B,D}$}
\begin{lem}\label{lem:smoothreal}
Let $B$, $D^\prime \in X(\R)$ such that the variety 
$Y_{B,D^\prime} \subset X \subset \PP^{n+1}$
given by \eqref{eqn:intersection} has a smooth real point $C^\prime$.
Let $V \subseteq X(\R)$ be an open neighbourhood of $C^\prime$.
Then there is an open neighbourhood $U \subseteq X(\R)$ of $D^\prime$,
such that for every $D \in U$, the variety $Y_{B,D}$
has a smooth real point $C \in V$.
\end{lem}
\begin{proof}
We may suppose that $B$, $C^\prime$, $D^\prime$ are contained in the affine patch $x_0=1$.
Let $G_1$, $G_2$, $G_3$ be the three polynomials defining $Y_{B,D^\prime}$
in \eqref{eqn:intersection} and let $g_1$, $g_2$, $g_3 \in
\R[x_1,\dotsc,x_{n+1}]$ be their dehomogenizations by $x_0=1$.
Write $f_1$, $f_2$, $f_3$ for the corresponding polynomials
in $\R[x_1,\dotsc,x_{n+1}]$ defining $Y_{B,D} \cap \{x_0=1\}$ with
$D \in X(\R) \cap \{x_0=1\}$. Of course $f_1=g_1$, $f_3=g_3$, and moreover
\[
\lVert f_2 - g_2 \rVert \le \mu \cdot \lVert  D- D^\prime \rVert_\infty
\]
where $\mu>0$ is a constant and $\lVert \cdot \rVert_\infty$ denotes
the infinity norm in the affine patch $x_0=1$ (which we identify
with $\R^{n+1}$). Now $C^\prime \in \R^{n+1}$ is 
a common zero
for $g_1$, $g_2$, $g_3$ with $\nabla{g_1}(C^\prime)$, $\nabla{g_2}(C^\prime)$,
$\nabla{g_3}(C^\prime)$ linearly independent (as $C^\prime$ is now
a smooth point on the affine patch $Y_{B,D^\prime} \cap \{x_0=1\}$).
Let $\varepsilon>0$ be sufficiently small so that $\cB_\varepsilon(C^\prime) \cap X(\R)$
is contained in $V$. 
Applying Lemma~\ref{lem:nr}, we know that
if $\lVert D-D^\prime \rVert_\infty$ is sufficiently small
then there is a non-zero vector $C \in \cB_\varepsilon(C^\prime)$ that is
a common zero for $f_1$, $f_2$, $f_3$ with
$\nabla{f_1}(C)$, $\nabla{f_2}(C)$, $\nabla{f_3}(C)$
linearly independent. This completes the proof.
\end{proof}

\begin{lem}\label{lem:analytic1}
Suppose $n \ge 48$. Let $B \in X(\Q)$. Suppose $D^\prime \in X(\R)$
such that $Y_{B,D^\prime}$ has a smooth real point $C^\prime$. 
Then there is
a non-empty open $U \subseteq X(\R)$ such that if $D \in U \cap X(\Q)$,
then $Y_{B,D}(\Q) \ne \emptyset$.
\end{lem}
\begin{proof}
Let $U$ be as in Lemma~\ref{lem:smoothreal}. Then $Y_{B,D}$
is defined over $\Q$ and has a smooth real point for every
$D \in U \cap X(\Q)$. Now the lemma follows from 
Corollary~\ref{cor:BDH}. 
\end{proof}

\section{A Little Geometry}\label{sec:diffgeom}

\subsection{Lines on $X$}
The following is
well-known  (for a proof, see
\cite[Lemma 2.1]{Siksek}).
\begin{lem}\label{lem:line}
Let $\ell$ be a line contained in $X$ and $P \in \ell$. Then
$\ell \subset T_P{X}$.
\end{lem}

\subsection{The second fundamental form}
Let $P \in X$. Associated to $P$ is a quadratic form (well-defined
up to multiplication by a non-zero scalar)
known as the second fundamental form which we denote
by $\Pi_P{X}$, and which is defined
as the differential of the Gauss map (e.g.\ \cite{GH}, 
\cite[Chapter 17]{Harris}). For our purpose the following
explicit recipe
given in \cite[pages 369--370]{GH} is useful.
By carrying out a non-singular change of 
coordinates we may suppose $P$ is $(1:0:\dotsc:0)$,
and the tangent plane $T_P{X}$ to $X$ at $P$ 
given by $x_{n+1}=0$.
 Then $X$ has
the equation $F=0$ with 
\begin{equation}\label{eqn:local}
F=x_0^2 x_{n+1}+x_0 q(x_1,\dotsc,x_{n+1})+c(x_1,\dotsc,x_{n+1})
\end{equation}
where $q$ and $c$ are homogeneous of degree $2$ and $3$ respectively.
Write $z_1=x_1/x_0,\dotsc,z_{n+1}=x_{n+1}/x_0$.
We can take $z_1,\dotsc,z_n$ as local coordinates for
$X$ at $P$, and then $X$ is given by the local equation
\[
z_{n+1}=q^\prime(z_1,\dotsc,z_n)+(\text{higher order terms}).
\]
Here $q^\prime(z_1,\dotsc,z_n)=-q(z_1,\dotsc,z_n,0)$.
The second fundamental form $\Pi_P{X}$ is the quadratic
form $q^\prime(dz_1,\dotsc,dz_n)$ (up to scaling). We shall only be concerned
with the rank and signature of $\Pi_P{X}$, which are
precisely the rank and signature of $q(x_1,\dotsc,x_{n},0)$
and so we will take this as the second fundamental form.
We may therefore view it as the restriction of $q$ to $T_P{X}$.
The following follows easily from the above description and the implicit function
theorem.

\begin{lem}\label{lem:closepoints}
Suppose $\Pi_P{X}$ has full rank $n$. 
\begin{enumerate}
\item[(i)] If $\Pi_P{X}$ is definite then there is an open neighbourhood
$U \subseteq X(\R)$ such that $U \cap X_P(\R)=\{P\}$.
\item[(ii)] 
If $\Pi_P{X}$ is indefinite then for every 
open neighbourhood
$U\subseteq X(\R)$ of $P$ the intersection
contains a real manifold of dimension $n-1$.
\end{enumerate}
\end{lem}

\begin{lem}\label{lem:indefiniteopen}
There is
a non-empty subset $U_1 \subseteq X(\R)$, open in the
real topology, such that 
for $P \in U_1$ the second fundamental
form $\Pi_P{X}$ is indefinite of
full rank. 
\end{lem}
\begin{proof}
A theorem of Landsberg \cite[Theorem 6.1]{Landsberg}
asserts that at a general point on smooth hypersurface
of degree $\ge 2$, the second fundamental form has full rank.
Thus there is a Zariski 
open $\cU \subset X$ such that $\Pi_P{X}$ has full rank
for $P \in \cU$. 

A straightforward application of Bertini's Theorem
shows the existence
of a real $3$-dimensional linear subvariety $\Lambda \subset
\PP^{n+1}$ such that 
$X^\prime=\Lambda \cap X$ is a smooth real cubic surface.
A classical theorem of Schl\"{a}fli asserts that the number
of real lines on smooth real cubic surface is either
$3$, $7$, $15$ or $27$. Let $\ell \subset \Lambda \cap X$
be a real line. By \cite[Lemma 2.2]{Siksek} all but 
at most two points of $\ell(\R)$ are hyperbolic
for $X^\prime$. Let $Q \in \ell(\R)$ be a hyperbolic
point for $X^\prime$. 
The determinant of the second fundamental
form $\Pi_Q{X^\prime}$
is the Gaussian curvature of $X^\prime$ at $Q$, which is
negative. It follows that the binary quadratic form $\Pi_Q{X^\prime}$
is indefinite. Now $\Pi_Q{X^\prime}$ is the restriction
of $\Pi_Q{X}$ to $T_Q{X^\prime}$ and so $\Pi_Q{X}$
is indefinite. Thus there is a neighbourhood
$V \subseteq X(\R)$ of $Q$, open in the real topology,
such that $\Pi_P{X}$ is indefinite for $P \in V$.
Now $V$ is necessarily Zariski-dense in $X$.
Thus $V \cap \cU(\R)$ is non-empty (as well
as being open in
the real topology). The proof is complete
upon letting $U_1=V \cap \cU(\R)$.
\end{proof}

\subsection{The Hessian}
Given $P \in X$, the Hessian of $F$ evaluated at $P$ is given 
by the $(n+2) \times (n+2)$ matrix 
\[
H_F(P)=\left( \frac{\partial^2 F}{\partial x_i \partial x_j}(P) \right)_{i,j=0,\dots,n+1}.
\]
Of course the Hessian is well-defined up to multiplication by a non-zero scalar.
\begin{lem}\label{lem:Hessian}
Let $P \in X$ and suppose $\Pi_P(X)$ has full rank $n$. Then $H_F(P)$ has full rank
$n+2$.
\end{lem}
\begin{proof}
%
Starting from \eqref{eqn:local},
an easy computation shows that the determinant
of the Hessian at $P$ is (up to sign) the
determinant of $q(x_1,\dots,x_n,0)$.
\end{proof}

\subsection{Eckardt Points}
We call $P \in X$ an \emph{Eckardt point} if $X_P:=X \cap T_P{X}$
is a cone with vertex at $P$. 
Note that if $n=2$ and $P$ is an Eckardt point
then $X_P$ consists of three lines meeting at $P$;
in this case $\Pi_P{X}$ vanishes identically.

For a proof of the following
classical 
theorem 
see \cite[Section 2]{CS}. 
\begin{thm}\label{thm:Eckardt}
The set of Eckardt points on $X$ is finite.
\end{thm} 

\subsection{Components of a Real Cubic Hypersurface}
We summarize briefly some well-known facts about components
of real cubic hypersurfaces. Everything we need is 
actually contained in \cite[Section 4.3]{Viro}.
A smooth 
real cubic hypersurface has either one or two connected components.
If it has two connected components then one of these is two-sided,
and homeomorphic to $S^n$, and the other is one-sided
and homeomorphic to $\R P^{n}$. If a line intersects the two-sided
component then it intersects it in two points, and intersects the 
odd-sided component in one point.
\begin{lem}\label{lem:twocomponents}
Suppose $X(\R)$ has two connected components. Then
$\Pi_P{X}$
is definite of full rank for all $P$ belonging to the two-sided
component.
\end{lem}
\begin{proof}
Let $P$ be a point on the two-sided component.
Suppose $\Pi_P{X}$ is indefinite or not of full rank.
Then there is a real line $\ell \subset T_P{X}$ (along which $\Pi_P{X}$ vanishes)
that meets $X$ with multiplicity
$\ge 3$ at $P$. As this is impossible for points on the two-sided
component, we have a contradiction.
\end{proof}

\section{Mordell--Weil Generation: First Steps}
\label{sec:tangentplane}
\begin{prop}\label{prop:tangentplane}
Let $P \in X(\Q)$ be a non-Eckardt point. 
Then the set $X_P(\Q)$ (considered as a subset of $X(\Q)$)
is contained in $\spn(P)$.
\end{prop}

The following lemma follows from the definitions.
\begin{lem}\label{lem:noncone}
Let $P \in X(\Q)$ and let $Q \in X_P(\Q)$ be distinct from $P$.
Suppose the line
$\ell$ joining $P$ to $Q$ is not contained in $X$. 
Then $Q \in \spn(P)$.
\end{lem}

For the proof of Proposition~\ref{prop:tangentplane} it remains
to show that $Q \in \spn(P)$ in the case $\ell \subset X$.
For $n=2$ this is \cite[Lemma 3.2]{Siksek}, so we suppose for
the remainder of this section that $n \ge 3$. 
\begin{lem}\label{lem:hypsec}
Any hyperplane section of $X$ is absolutely irreducible.
\end{lem}
\begin{proof}
Let $L=0$ be a hyperplane such that $X \cap \{L=0\}$ is
absolutely reducible. Then we can write 
$F=LQ+L^\prime Q^\prime$ where $L$, $L^\prime$ are homogeneous
linear, and $Q$, $Q^\prime$ are homogeneous quadratic. As $n \ge 3$,
the variety
$L=L^\prime=Q=Q^\prime=0$
has a point $R \in \PP^{n+1}$. It follows that 
$R$ is a singular on $X$ giving a contradiction.
\end{proof}

\begin{lem}\label{lem:distinctplane}
Let $P \in X(\Q)$ be a non-Eckardt point. 
Let $Q \in X_P(\Q)$ with $T_Q{X} \ne T_P{X}$.
Then $Q \in \spn(P)$.
\end{lem}
\begin{proof}
Let $\cW \subseteq X_P$ be the subvariety consisting of
lines through $P$ contained in $X_P$.
As $P$ is a non-Eckardt point, $\cW$ is a proper
subvariety. Moreover, by Lemma~\ref{lem:hypsec},
the tangent plane section $X_P$
is irreducible, and so $\dim(\cW)< \dim(X_P)$.
Let $\cU=X_P-\cW$ which is Zariski dense in $X_P$.

Let $\cV:=X_P \setminus (X_P \cap T_Q)$.
As $T_Q \ne T_P$, this is a dense open subset of $X_P$.
Let $\iota \; : \; \cV \rightarrow \cV$ be the involution
given as follows. If $R \in \cV$ we join $R$ to $Q$ by the line $\ell_{R,Q}$
and we let $\iota(R)$ be the third point of intersection of this line with $X$.
We note that $\ell_{R,Q} \not \subset X$, since otherwise it will be contained
in $T_Q{X}$ by Lemma~\ref{lem:line}. Now $(\cV \cap \cU) \cap \iota (\cV \cap \cU)$ 
is a Zariski dense subset of the rational
variety $X_P$. This dense subset must contain a rational point 
$R$. Then $R$, $\iota(R) \notin \cW$ and so $R$, $\iota(R) \in \spn(P)$
by Lemma~\ref{lem:noncone}. Finally the line joining $R$ with $\iota(R)$
passes through $Q$ and is not contained in $X$. Thus $Q \in \spn(R)$.
\end{proof}
 
\begin{proof}[Proof of Proposition~\ref{prop:tangentplane}]
Let $Q \in X_P(\Q)$. We would like to show that $Q \in \spn(P)$.
Thanks to Lemmas~\ref{lem:noncone} and~\ref{lem:distinctplane},
we may suppose there is a $\Q$-line $\ell \subset X$ containing $P$, $Q$,
and $T_Q {X}= T_P{X}$. Now the line $\ell$ contains at most finitely
many Eckardt points by Theorem~\ref{thm:Eckardt}. Moreover, the Gauss map
on a smooth hypersurface
has finite fibres \cite[Lecture 15]{Harris}. 
Thus there is a non-Eckardt $R \in \ell(\Q)$ with $T_R{X} \ne T_P{X}$.
It follows that $R \in \spn(P)$.
Moreover, $Q \in \ell \subset T_R$ by Lemma~\ref{lem:line} and so
$Q \in \spn(R)$ (again by Lemma~\ref{lem:distinctplane}). This completes the proof.
\end{proof}

\begin{lem}\label{lem:analytic2}
Suppose $n \ge 48$. Let $B \in X(\Q)$ so
that $X_B$ does not contain points that are Eckardt for $X$. 
Suppose $D^\prime \in X(\R)$
such that $Y_{B,D^\prime}$ has a smooth real point $C^\prime$. 
Then there is
a non-empty open $U \subseteq X(\R)$ such that $U \cap X(\Q) \subseteq \spn(B)$.
\end{lem}
\begin{proof}
Take $U$ to be as in Lemma~\ref{lem:analytic1}.
Let $D \in U \cap X(\Q)$. By the conclusion of Lemma~\ref{lem:analytic1}
we see that $Y_{B,D}$ has a rational point $C$. From the equations
defining $Y_{B,D}$ in \eqref{eqn:intersection} we have that
$C \in X_B(\Q)$ and $D \in X_C(\Q)$. Moreover,  
neither $B$ nor $C$ (both contained in $X_B$)
are Eckardt points. Applying Proposition~\ref{prop:tangentplane},
we have $C \in \spn(B)$ and $D \in \spn(C)$ completing the proof.
\end{proof}

\section{A Smoothness Criterion}

\begin{lem}\label{lem:smooth}
Let $B \in X$, 
Let $C^\prime \in X_B$ and $D^\prime \in X_{C^\prime}$. Suppose 
\begin{enumerate}
\item[(i)] $T_{C^\prime}{X} \ne T_B{X}$; 
\item[(ii)] $H_F(C^\prime)$ has full rank, where $H_F$ is 
the Hessian matrix;
\item[(iii)] $D^\prime$ does not belong to the line
\begin{equation}\label{eqn:line}
\{ \; (\lambda \nabla{F}(B)+\mu \nabla{F}(C^\prime))\cdot
H_F(C^\prime)^{-1} 
\quad : \quad (\lambda : \mu) \in \PP^1 \; \} \, .
\end{equation}
\end{enumerate}
Then $C^\prime$ is a smooth point on the variety $Y_{B,D^\prime} \subset \PP^{n+1}$
given by \eqref{eqn:intersection}.
\end{lem}
\begin{proof}
As $C^\prime \in X_B$ and $D^\prime \in X_{C^\prime}$ we see that 
$C^\prime \in Y_{B,D^\prime}$. We need to show that $C^\prime$ is a smooth point on
$Y_{B,D^\prime}$. Write
\[
f(x_0,\dotsc,x_{n+1})=\nabla{F}(x_0,\dotsc,x_{n+1}) \cdot D^\prime,\qquad
g=\nabla{F}(B) \cdot (x_0,\dotsc,x_{n+1}).
\] 
To show that $C^\prime$ is smooth on $Y_{B,D^\prime}$ it is enough to show that
$\nabla{F}(C^\prime)$, $\nabla{f}(C^\prime)$ and $\nabla{g}(C^\prime)$ are linearly
independent. A straightforward computation shows that
\[
\nabla{f}(C^\prime)=D^\prime \cdot H_F(C^\prime), \qquad \nabla{g}(C^\prime)=\nabla{F}(B).
\]
Suppose 
\[
\varepsilon D^\prime \cdot H_F(C^\prime)+\lambda \nabla{F}(B)+\mu \nabla{F}(C^\prime)=0.
\]
By assumptions (ii) and (iii) we see that $\varepsilon=0$.
However, $\nabla{F}(B)$ and $\nabla{F}(C^\prime)$
are linearly independent by assumption (i), and so $\lambda=\mu=0$. 
\end{proof}

\section{Proof of Theorem~\ref{thm:main}}
In this section $n \ge 48$.
\begin{lem}\label{lem:initialpatch}
There is $A \in X(\Q)$ and a non-empty
open $U \subseteq X(\R)$ such that:
\begin{enumerate}
\item[(i)] $U \cap X(\Q) \subseteq \spn(A)$;
\item[(ii)] $\spn(A)$ contains at least one point in every
connected component of $X(\R)$.
\end{enumerate}
\end{lem}
\begin{proof}
Suppose first that $X(\R)$ is connected.
Let $U_1 \subseteq X(\R)$ be the non-empty open subset whose
existence is guaranteed by Lemma~\ref{lem:indefiniteopen}: for
every $P \in U_1$, the second fundamental form $\Pi_P{X}$ is
indefinite of full rank. It follows from Theorem~\ref{thm:Eckardt}
that the set of points $P$ with $X_P$ containing an Eckardt point 
is a proper subset of $X$ that is closed in the Zariski
topology. Thus we may replace $U_1$
by a non-empty open set $U_2 \subseteq U_1$ such that for every $P \in U_2$,
the subvariety $X_P$ does not points that are Eckardt for $X$.
Fix
$B \in U_2 \cap X(\Q)$ whose existence is guaranteed by Theorem~\ref{thm:Skinner}. 
The hypersurface $X$ is smooth of degree $3$,
and so
the Gauss map $X \rightarrow X^*$
has finite fibres \cite[Lecture 15]{Harris}. 
We can therefore take an open neighbourhood $U_3 \subseteq U_2$
of $B$ such that for all $C^\prime \in U_3$ with $C^\prime \ne B$, 
we have $T_B{X} \ne T_{C^\prime}{X}$.
By Lemma~\ref{lem:closepoints}, the intersection
$U_3 \cap X_B(\R)$ contains 
a real manifold of dimension $n-1$; choose $C^\prime \in U_3 \cap X_B(\R)$
with $C^\prime \ne B$.
As the second fundamental form has full
rank on $U_3$, we see from Lemma~\ref{lem:Hessian} that
$H_F(C^\prime)$ is of full rank $n+2$.
 Now again by Lemma~\ref{lem:closepoints}, the intersection
$U_3 \cap X_{C^\prime}(\R)$ contains a manifold of real dimension $n-1$,
and so we can find $D^\prime \in U_3 \cap X_{C^\prime}(\R)$ that avoids the line
\eqref{eqn:line}. The points $B$, $C^\prime$, $D^\prime$
satisfy the conditions of Lemma~\ref{lem:smooth}.
Thus $C^\prime$ is a smooth point on $Y_{B,D^\prime}$.
By Lemma~\ref{lem:analytic2}, there is a non-empty open $U$ such that $U \cap X(\Q) \subseteq \spn(B)$.
We simply take $A=B$, and the proof is complete in the case
when $X(\R)$ is connected.


Now suppose $X(\R)$ has two connected components. Let $U_2 \subseteq X(\R)$
be as above. From Lemma~\ref{lem:twocomponents} we know that $U_2$ is contained
in the one-sided component. Let $B^\prime \in U_2$. Let $\ell^\prime$ be a real
line passing through $B^\prime$ and tangent to the two-sided component
at a point $A^\prime$. 
By Corollary~\ref{cor:tangentbundle}, there is a point $A \in X(\Q)$
belonging to the two-sided component and a line $\ell$ defined over
$\Q$ such that $\ell \cdot X=2A+B$ where $B \in U_2 \cap X(\Q)$.
Now $B \in \spn(A)$
and $\spn(A)$ contains points belonging to both components of $X(\R)$.
From the above argument there is a non-empty open $U \subseteq X(\R)$ such that
$U \cap X(\Q) \subseteq \spn(B) \subseteq \spn(A)$.
\end{proof}


\begin{lem}\label{lem:patchenlarge}
Let $A \in X(\Q)$ be as in Lemma~\ref{lem:initialpatch}.
Then there is an open $W \subseteq X(\R)$ such that
$W \cap X(\Q)=\spn(A)$.
\end{lem}
\begin{proof}
Let $U$ be as in Lemma~\ref{lem:initialpatch}.
We may suppose $\spn(A) \not \subset U$, otherwise
we simply take $W=U$ and
there is nothing to prove. 
Let $P \in \spn(A)$ that does not belong to $U$. 
By Theorem~\ref{thm:Skinner}, 
there is some 
$P^\prime \in U \cap X(\Q)$ such that $P^\prime \notin T_{P}X$.
Let $\ell$ be the line joining $P$ to $P^\prime$. 
The line 
$\ell$ is not contained in $T_P{X}$ and so,
by Lemma~\ref{lem:line}, not contained in $X$.
Let $P^{\prime\prime} \in X(\Q)$ be the third point of intersection of $\ell$
with $X$. Since $P \in \spn(A)$ and $P^\prime \in U \cap X(\Q) \subseteq \spn(A)$, we have
$P^{\prime\prime} \in \spn(A)$.
Observe that $\ell$ is not contained
in the tangent plane of $P^{\prime\prime}$ (for otherwise
$\ell$ would be contained in $X$).
Now there is some non-empty open $U^\prime \subset U$
containing $P^\prime$ that is disjoint from the tangent
plane of $P^{\prime\prime}$. For a point $R \in U^\prime$,
let $\varphi(R)$ denote the third point of intersection
of the (real) line joining $R$ to $P^{\prime\prime}$. Then the map
$\varphi \; : \; U^\prime \rightarrow X(\R)$ is continuous
and injective. By the Invariance of Domain Theorem 
\cite[Corollary IV.19.9]{Bredon},
the image $\varphi(U^\prime)$ is open. 
We shall let $W_P=\varphi(U^\prime)$.
Clearly $P \in W_P$ and $W_P \cap X(\Q)
\subseteq
\spn(A)$. The lemma follows on taking 
\[
W=U \cup \bigcup_{P \in \spn(A) \setminus U} W_P.
\]
\end{proof}

\begin{lem}\label{lem:secantclosed}
Let $W$ be as in Lemma~\ref{lem:patchenlarge}, and write
$\overline{W}$ for its closure. 
Then $\overline{W}$ is closed under 
secant operations: if $P$, $Q \in \overline{W}$
are distinct, and if the line $\ell$ joining them is not contained
in $X$, then $R \in \overline{W}$ where $\ell \cdot X=P+Q+R$.
\end{lem}
\begin{proof}
By Theorem~\ref{thm:Skinner} there
exist $\{P_k\}$, $\{Q_k\} \subset W \cap X(\Q)$, with $P_k \ne Q_k$, that converge respectively to $P$, $Q$.
Write
$F \subset \mathbb{G}(n+1,1)$ for the Fano scheme
of lines on $X$. Then the real points of $F$ are closed
in $\mathbb{G}(n+1,1)(\R)$.
As $\ell \notin F(\R)$, we see for large enough $k$
that the line $\ell_k/\Q$ joining $P_k$, $Q_k$ is not contained
in $X$. Let $\ell_k \cdot X=P_k+Q_k+R_k$.
Then $\{R_k\}$ converges to $R$. Moreover, $P_k$, $Q_k \in 
W \cap X(\Q) \subseteq \spn(A)$.
Hence $R_k \in \spn(A) \subset W$ and so $R \in \overline{W}$. 
\end{proof}

\begin{lem}\label{lem:components}
Let $A$, $W$ be as above.
Then $\overline{W}=X(\R)$.
\end{lem}
\begin{proof}
We claim that $\overline{W}$ is open. From that it follows that
$\overline{W}$ is a union of connected components of $X(\R)$.
As $\spn(A) \subset \overline{W}$ contains points from every component, the
lemma follows from the claim.

To prove the claim we mimic the argument in the proof of Lemma~\ref{lem:patchenlarge}.
Let $P\in \overline{W}$. 
Let 
$P^\prime \in W$ such that $P^\prime \notin T_{P}X$,
and let $\ell$ be the line joining $P$ to $P^\prime$.  As $\overline{W}$
is closed under secant operations, 
$P^{\prime\prime} \in \overline{W}$ where $\ell \cdot X=P+P^\prime+P^{\prime\prime}$.
Now there is some non-empty open $W^\prime \subset W$
containing $P^\prime$ that is disjoint from the tangent
plane of $P^{\prime\prime}$. For a point $R \in W^\prime$,
let $\varphi(R)$ denote the third point of intersection
of the (real) line joining $R$ to $P^{\prime\prime}$. Then the map
$\varphi \; : \; W^\prime \rightarrow X(\R)$ is continuous
and injective, and thus 
the image $\varphi(W^\prime)$ is open. Clearly $\varphi(W^\prime)$
contains $P$ and is contained in $\overline{W}$ (as the latter
is closed under secant operations).
\end{proof}

\begin{proof}[Proof of Theorem~\ref{thm:main}]
Let $A$, $W$ be as above. In particular, $\spn(A)=W \cap X(\Q)$
and $\overline{W}=X(\R)$.
We write
$\partial{W}=X(\R)\setminus W$. We note that $\partial{W}$
is the complement of an open dense set, and therefore nowhere
dense. 

We want to show that $X(\Q)=\spn(A)$. Let
$P \in X(\Q)$.  
Then there is a Zariski open $\cU \subset X$
and an involution $\iota \; : \; \cU \rightarrow \cU$ that
sends $R \in \cU$ to the third point of the line joining $R$ to $P$.
Choose $R \in \cU(\R) \cap X(\Q)$ such that $R \notin \partial{W} \cup \iota (\partial{W})$.
Then $R$, $\iota(R) \in \spn{A}$ and so $P \in \spn(A)$.
\end{proof}

\end{document}